\documentclass[12pt]{amsart}
\newif\ifpdf\ifx\pdfoutput\undefined\pdffalse\else\pdfoutput=1\pdftrue\fi

\usepackage{graphicx,amsmath,amsfonts,latexsym,amssymb,amsthm,amscd}
\usepackage{latexsym}


\newcommand{\field}[1]{\mathbb{#1}}
\newcommand{\rz}{\field{R}}
\newcommand{\cz}{\field{C}}
\newcommand{\zz}{\field{Z}}
\newcommand{\nz}{\field{N}}

\newcommand{\I}{\mathrm{i}}

\newcommand{\Hi}{\mathcal{H}}

\newcommand{\spec}{\mathrm{spec}}

\newcommand{\lieg}{\mathfrak{g}}
\newcommand{\liek}{\mathfrak{k}}
\newcommand{\liea}{\mathfrak{a}}
\newcommand{\lieh}{\mathfrak{h}}
\newcommand{\liep}{\mathfrak{p}}
\newcommand{\liem}{\mathfrak{m}}
\newcommand{\lien}{\mathfrak{n}}
\newcommand{\liegc}{\mathfrak{g}_\cz}

\newcommand{\Ad}{\mathrm{Ad}}
\newcommand{\killing}{\kappa}
\newcommand{\hc}{\mathbf{c}}

 \newtheorem{theorem}{Theorem}[section]

 \newtheorem{pro}[theorem]{Proposition}

\title{Analytic Continuation of Resolvent Kernels on noncompact Symmetric Spaces}

\author{Alexander Strohmaier}

\begin{document}

\maketitle

\begin{abstract}
  \noindent
  Let $X=G/K$ be a symmetric space of noncompact type and let $\Delta$ be the
  Laplacian associated with a $G$-invariant metric on $X$.  We show that the
  resolvent kernel of $\Delta$ admits a holomorphic extension to a Riemann
  surface depending on the rank of the symmetric space. This Riemann
  surface is a branched cover of the complex plane with a certain part of the
  real axis removed.  It has a branching point at the bottom of the spectrum
  of $\Delta$.  It is further shown that this branching point is quadratic if
  the rank of $X$ is odd, and is logarithmic otherwise. In case $G$ has only
  one conjugacy class of Cartan subalgebras the resolvent kernel extends to a
  holomorphic function on a branched cover of $\cz$ with the only branching
  point being the bottom of the spectrum.
\end{abstract}



\section{Introduction}

If $(M,g)$ is a complete Riemannian manifold the metric Laplace operator
$\Delta_g$ is a selfadjoint operator in the Hilbert space $L^2(M)$.  Its
resolvent
\begin{gather}
  R_z(\Delta_g)=(\Delta_g -z)^{-1}
\end{gather}
is a holomorphic function on $\cz \backslash \textrm{spec}(\Delta_g)$ with
values in the bounded operators in $L^2(M)$.  Whereas $\cz \backslash
\textrm{spec}(\Delta_g)$ is the largest domain where all matrix elements of
the resolvent
\begin{gather*}
  \langle \phi, R_z(\Delta_g) \psi \rangle \quad \phi,\psi \in L^2(M)
\end{gather*}
are holomorphic functions it may happen that for all $\psi$ and $\phi$ in some
dense subset of $L^2(M)$ those matrix elements can be continued to holomorphic
functions on a Riemann surface extending the resolvent set $\cz \backslash
\textrm{spec}(\Delta_g)$.

A lot of examples where such a situation occurs are known. If $(M,g)$ is the
Euclidean space of dimension $n$, then the functions
\begin{gather*}
  \langle \phi, R_z(\Delta_g) \psi \rangle \quad \phi,\psi \in
  L^2(\rz^n,e^{+\epsilon x^2} dx), \quad \epsilon >0
\end{gather*}
extend to holomorphic functions on the concrete Riemann surface associated
with either the function $\sqrt{z}$ or $\log(z)$, depending on whether the
dimension $n$ is odd or even (see e.g. \cite{MR96k:35129}).  
If a compactly supported potential is added,
the same holds, except that in this case the continuation is meromorphic.
Poles of these matrix elements are commonly referred to as scattering poles or
resonances.  The theory that provides analytic continuation in the case of the
Laplacian in $\rz^n$ with potential is well developed and plays a crucial role
for example in proving the absence of singular continuous spectrum (see e.g.
\cite{CFKS} and references therein).  Also in
order to make sense of the notion of resonances one needs to proof the
existence of a meromorphic continuation across the spectrum.  From the results
known for the Euclidean space one may obtain further examples employing
perturbation techniques (e.g. \cite{MR96j:35183}).

Here we will be mainly interested in a continuation of the resolvent kernel,
i.e. we view $R_z(\Delta_g)$ as a function with values in $\mathcal{D}'(M
\times M)$. Hence, we look at the matrix elements $\langle \phi, R_z(\Delta_g)
\psi \rangle$ where $\phi,\psi \in C^\infty_0(M)$.  

Examples of manifolds where a meromorphic continuation of the resolvent kernel
is known to exist include manifolds with cylindrical ends 
(\cite{MR96k:35129,MR98f:11053,MR96g:58180}), 
certain strictly pseudoconvex domains (\cite{MR92i:32016}), 
manifolds with asymptotically constant negative curvature
(\cite{MR89c:58133,MR96h:58172}), and noncompact symmetric spaces of rank $1$ 
(\cite{MR93i:58160,CarrPed:03}). Concerning this question little was known
however about higher rank symmetric spaces of noncompact type. 
Only recently the method of complex scaling
was used by Mazzeo and Vasy (\cite{MazVas:02,MazVas:03}) to show that 
the resolvent kernel has a meromorphic extension across the spectrum to a 
region which is obtained by rotating the spectrum up to the angle $\pi/2$ and
by adding further branching points.
It was conjectured however in their papers that this continuation has no poles
and that the surface can be chosen without branching points.
A special case which can be considered as well understood is that
 of a symmetric space $X=G/K$, where $G$ is a complex semisimple Lie group. 
In this case there is an explicit
formula for the heat kernel (see \cite{MR39:360}). By Laplace transform
one can obtain an explicit expression for the resolvent kernel in terms of 
the modified Bessel function of the second kind. 
In this way one finds that the resolvent kernel admits a holomorphic extension
to a branched cover of the complex plane.

We will show that in the case of a noncompact Riemannian symmetric space
$X=G/K$ of arbitrary rank the resolvent kernel admits a holomorphic extension
to a Riemann surface depending on the rank of the symmetric space. This
Riemann surface is a branched cover of the complex plane with a part of the
real axis removed. A branching point occurs at the bottom of the spectrum.
This branching point turns out to be quadratic if the rank of the space is odd
and logarithmic otherwise.  This Riemann surface contains the region
obtained by rotating the spectrum up to the angle $\pi$, and our result
implies the absence of poles in that region.
In case the rank of the space is one we reproduce the known
results, i.e.  the resolvent kernel extends to a meromorphic function on the
branched cover of $\cz$ associated with the function $\sqrt{z-\mu}$,
where $\mu$ is the bottom of the spectrum.  We get even
stronger results in case the group $G$ has only one conjugacy class of Cartan
subalgebras. In this case the resolvent kernel behaves almost like in
Euclidean space, i.e. there exists a holomorphic continuation to a Riemann
surface which is a branched cover of $\cz$ with the only branching point
being the bottom of the spectrum.
The branching point is quadratic or logarithmic depending on
whether the rank of the space is odd or even.  This special case includes the
real hyperbolic spaces of odd dimension, the spaces $SU^*(2n)/Sp(n)$,
$E_{6(-26)}/F_4$ as well as the case $G/K$ when $G$ is a complex group.

Our method relies on the Fourier transform on symmetric spaces of noncompact
type and on the meromorphic continuation of the Harish-Chandra $\hc$-function.

\section{Notations and Background material}

\subsection{Symmetric spaces of noncompact type}
Suppose that $X=G/K$ is a symmetric space of noncompact type, that is $G$ is a
real connected noncompact semisimple Lie group with finite center and $K
\subset G$ is a maximal compact subgroup.  Let $\lieg$ be the Lie algebra of $G$ and let
$\liek \subset \lieg$ be the Lie algebra of $K$.  We denote by $\killing(\cdot,\cdot)$
the Killing form on $\lieg$, i.e.
\begin{gather}
\killing(X,Y)=\mathrm{Tr}(\Ad(X) \circ \Ad(Y)).
\end{gather}
Then there exists a Cartan involution $\theta: \lieg \to \lieg$ with fixed
point algebra $\liek$, i.e. $\theta$ is an involutive automorphism such that
the bilinear form $\langle X,Y \rangle:=-\killing(X,\theta Y)$ is positive definite
and such that the $+1$ eigenspace of $\theta$ coincides with $\liek$.  Hence,
the decomposition of $\lieg$ into $+1$ and $-1$ eigenspaces of $\theta$ reads
\begin{gather}
\lieg = \liek \oplus \liep,
\end{gather}
where $\liep \subset \lieg$ is a linear subspace.  If $\liegc$ is the
complexification of $\lieg$ then
\begin{gather}
\lieg_c=\liek + \I \liep \subset \liegc
\end{gather}
is a compact real form of $\liegc$.  Now let $\liea$ be a maximal abelian
subspace of $\liep$ and let $\liem$ be the centralizer of $\liea$ in $\liek$.
If $\mathfrak{b}$ is a maximal abelian subalgebra of $\liem$ then $\lieh =
\mathfrak{b} + \liea$ is a $\theta$-stable Cartan subalgebra of $\lieg$. The
set of restricted roots $\Delta(\lieg,\liea)$ w.r.t. $\liea$ is the set of
nonzero linear functionals $\alpha \in \liea^*$ such that
\begin{gather}
[x,y]=\alpha(x) y, \quad \forall x \in \liea
\end{gather}
for some nonzero element $y \in \lieg$. The multiplicity $m_\alpha$ of a
restricted root is the dimension of the vector space $\{y \in \lieg;\;
[x,y]=\alpha(x) y, \quad \forall x \in \liea\}$.  This gives rise to a root
space decomposition of $\lieg$:
\begin{gather}
\lieg = \liea \oplus \liem \oplus \sum_{\alpha \in \Delta(\lieg,\liea)}
\lieg_\alpha,
\end{gather}
where $\lieg_\alpha$ are the root subspaces, i.e.
\begin{gather}
\lieg_\alpha=\{y \in \lieg;\; [x,y]= \alpha(x) y, \quad \forall x \in
\liea\}.
\end{gather}
Note that each restricted root $\Delta(\lieg,\liea)$ coincides with the restriction
of a root $\alpha \in \Delta(\lieg,\lieh)$ to $\liea$. We may choose now a subsystem of positive roots
$\Delta^+(\lieg,\lieh)$ and accordingly a subsystem of positive restricted
roots $\Delta^+(\lieg,\liea)$.  As usual we use the notation
\begin{gather}
  \rho=\frac{1}{2} \sum_{\alpha \in \Delta^+(\lieg,\liea)} m_\alpha \alpha
\end{gather}
for half of the sum of positive roots counting multiplicities.  Let $\lien
\subset \lieg$ be the nilpotent Lie-subalgebra
\begin{gather}
\lien=\sum_{\alpha \in \Delta^+(\lieg,\liea)} g_\alpha.
\end{gather}
Then the Iwasawa decomposition of $\lieg$ reads as follows:
\begin{gather}
\lieg = \liek \oplus \liea \oplus \lien.
\end{gather}

Now let $A$ and $N$ be the analytic subgroups of $G$ with Lie algebras $\liea$
and $\lien$ respectively. Then
\begin{gather}
N \times A \times K \to G, \; (n,a,k) \to nak
\end{gather}
is a diffeomorphism from $N \times A \times K$ to $G$.  This is the
classical Iwasawa decomposition of $G$.  For $g \in G$ we denote by $A(g) \in
\liea$ the unique element such that $g$ can be expressed as
\begin{gather}
  g = n \exp(A(g)) k
\end{gather}
with $n \in N$ and $k \in K$.
Let $M$ be the centralizer of $A$ in $K$ and denote by $B$ be the compact homogeneous
space $K/M$.  Then $A(gK,kM):=A(k^{-1} g)$ defines a smooth function on $X
\times B$ with values in $\liea$.  
If $M'$ is the normalizer of $A$ in $K$ then the restricted Weyl group is
defined by $W:=M'/M$. The restricted Weyl group acts on $\liea$ by 
$kM a = \mathrm{Ad}(k) a$ and by duality it also acts on the dual $\liea^*$ of
$\liea$. This representation of $W$ on $\liea$ is injective and $W$ can be
identified in this way with a group of reflections in $\liea$.
As usual let $\liea_+^*$ be the positive Weyl chamber 
$\{\lambda \in \liea^*;\; \lambda(\alpha)>0 \quad \forall \alpha \in \Delta^+(\lieg,\liea)\}$.
Then $\liea_+^*$ is a fundamental domain for the action of $W$ on $\liea^*$.

\subsection{The Helgason transform and the Paley-Wiener theorem}

Let $X=G/K$ be a symmetric space of noncompact type. We use the notations and
conventions from above.  The generalized Fourier transform (Helgason
transform) of a function $f \in C^\infty_0(X)$ is the function $\hat f :
\liea^* \times B \to \cz$ defined by
\begin{gather} 
  \hat f(\lambda,b)=\int_X  f(x) e^{(-\I \lambda + \rho)(A(x,b))} dx,
\end{gather}
where we integrate with respect to some invariant measure on $X$ which we
choose to be normalized as in \cite{MR41:8587}
\footnote{the normalization of the measure is not important for our
  considerations}.
For suitable normalizations of the Euclidean measure on $\liea^*$ and the invariant
measure $db$ on $B$ the inverse Fourier transform (see \cite{MR41:8587}) is given by
\begin{gather} 
  f(x)= w^{-1} \int_{\liea^* \times B}
  \hat f(\lambda,b) e^{(+\I \lambda + \rho)(A(x,b))} \frac{d
  \lambda}{|\hc(\lambda)|^2} db,
\end{gather}
where $\hc(\lambda)$ is the Harish-Chandra $\hc$-function, and $w$ is the order of
the restricted Weyl group. Moreover, for $f,g \in C^\infty_0(X)$ we have
\begin{gather}
 \int_X  \overline f(x) g(x) dx = w^{-1} \int_{\liea^* \times B}
 \overline{\hat f}(\lambda,b) \hat
 g(\lambda,b) \frac{d \lambda}{|\hc(\lambda)|^2} db .
\end{gather}
The Fourier transform $\mathcal{F}$ extends to an
isometry
\begin{gather*}
  L^2(X) \to L^2(\liea^*_+ \times B,\frac{d
    \lambda}{|\hc(\lambda)|^2} db).
\end{gather*}
For $f \in C^\infty_0(X)$ the Fourier transform $\hat f$ extends to a function on
$\liea^*_\cz \times B$ which is entire in the first variable.  Let
$\Hi(a^*_\cz)$ be the space of holomorphic functions of uniform exponential
type on $a^*_\cz \times B$,
i.e. the space functions $\phi$ on $a^*_\cz \times B$, entire in the
first variable, such that
there exists a constant $R>0$ with
\begin{gather}
  |\phi(\lambda,b)| \leq C_N (1+|\lambda|_\killing)^{-N} e^{R |\Im(\lambda)|_\killing}
\end{gather}
for all $N \in \nz_0$ with some constants $C_N>0$. Here $| \cdot|_\killing$ denotes
the norm induced by the Killing form. 
The Paley-Wiener theorem for symmetric spaces of noncompact type (see
\cite{MR51:3804}) is precisely the statement that the image of $C^\infty_0(X)$
under the Fourier transform coincides with the space of functions
$f \in \Hi(a^*_\cz)$ which satisfy
\begin{gather}
 \int_B f(s\lambda,b) e^{(+\I s\lambda + \rho)(A(x,b))} db= 
 \int_B f(\lambda,b) e^{(+\I \lambda + \rho)(A(x,b))} db \quad \forall s
 \in W, x \in X, \lambda \in \liea^*_\cz.
\end{gather}

The $\hc$-Function satisfies $|\hc(\lambda)|^2=\hc(\lambda) \hc(-\lambda)$ and is
known to extend to a a meromorphic function on $\liea^*_\cz$.  This follows
from the product formula of Gindikin and Karpelevic (\cite{MR27:240}) which
leads to the following explicit expression for the $\hc$-function.
 \begin{gather} \label{GKFormula}
   \hc(\lambda)=\hc_0 \prod_{\alpha \in \Delta(\lieg,\liea)_0^+} \frac{2^{-\langle
       \I \lambda, \alpha_0 \rangle} \Gamma(\langle \I \lambda, \alpha_0
     \rangle)} {\Gamma(\frac{1}{2}(\frac{1}{2} m_\alpha + 1 + \langle \I
     \lambda, \alpha_0 \rangle)) \Gamma(\frac{1}{2}(\frac{1}{2} m_\alpha +
     m_{2 \alpha} + \langle \I \lambda, \alpha_0 \rangle))},
 \end{gather}
 where $\alpha_0 = \frac{\alpha}{\langle \alpha,\alpha \rangle}$ and the
 constant $\hc_0$ is determined by $\hc(-\I \rho)=1$. The product is over the set
 $\Delta(\lieg,\liea)_0^+$ of indecomposable positive restricted roots, i.e.
 \begin{gather}
   \Delta(\lieg,\liea)_0:=\{\alpha \in \Delta(\lieg,\liea)^+;\; \frac{1}{2}
   \alpha \notin \Delta(\lieg,\liea)\}.
 \end{gather}
 The scalar product $\langle \cdot,\cdot\rangle$ is the positive definite
 scalar product on $\liea^*$ induced by the Killing form.
 
\subsection{The metric Laplacian}

As above let $X=G/K$ be a symmetric space of noncompact type and let
$\lieg=\liek  \oplus \liep$ the corresponding Cartan decomposition.  The
tangent space $T_{eK} X$ at the point $eK$ in $X$ is canonically identified
with $\liep$ and invariant metrics on $X$ are therefore in $1-1$
correspondence with $Ad(K)$-invariant positive definite scalar products on
$\liep$. Suppose now that $g$ is an invariant metric on $X$ and denote by
$(\cdot,\cdot)$ the corresponding scalar product on $\liep$.  This induces a
scalar product on $\liea$ and on its dual $\liea^*$.  We denote the dual
scalar product on $\liea^*$ again by $(\cdot,\cdot)$ \footnote{in contrast to
  the scalar product induced by the Killing form which we denoted by $\langle
  \cdot,\cdot\rangle$}.  The Laplace operator $\Delta_g: C^\infty_0(X) \to
L^2(X)$ associated with the metric is a $G$-invariant differential operator
and hence, its spectral decomposition is strongly related to the Fourier
transform.  For $\lambda \in \liea^*$ and $b \in B$ we define the function
$\phi_{\lambda,b}$ on $X$ by $\phi_{\lambda,b}(x)=e^{(-\I \lambda +
\rho)(A(x,b))}$ and obtain (see \cite{MR20:925} and \cite{MR51:3804} p. 458, 462)
\footnote{this formula appears in the literature for the case that the metric is induced
by the Killing form. On each irreducible factor of the symmetric space any
invariant metric is proportional to the one induced by the Killing form. Taking this into
account it is easy to see that the formula holds in the general situation.}
\begin{gather}
  \Delta_g \phi_{\lambda,b} = ( |\lambda|^2 + |\rho|^2 ) \phi_{\lambda,b},
\end{gather}
where the norm $|\lambda|=( \lambda ,\lambda )^{1/2}$ is with respect to the
scalar product induced by the metric $g$.  As a consequence if $\mathcal{F}:L^2(X) \to
L^2(\liea^*_+ \times B,\frac{d \lambda}{|\hc(\lambda)|^2} db)$ is the
Fourier transform the operator $\mathcal{F} \Delta_g \mathcal{F}^{-1}$ is a
multiplication operator.  More precisely, if $m_{|\lambda|^2 + |\rho|^2}$ is
the selfadjoint operator in $L^2(\liea^*_+ \times B, \frac{d
  \lambda}{|\hc(\lambda)|^2} db)$ defined by
\begin{gather} \label{multi}
  (m_{|\lambda|^2 + |\rho|^2}f)(\lambda,k)= ( |\lambda|^2 + |\rho|^2 )
  f(\lambda,k).
\end{gather}
with maximal domain, the operator $\mathcal{F} \Delta_g \mathcal{F}^{-1}$
coincides with its restriction to the set $\mathcal{F}(C^\infty_0(X))$.  Therefore,
$m_{|\lambda|^2 + |\rho|^2}$ coincides with $\mathcal{F} \overline{\Delta_g}
\mathcal{F}^{-1}$, where $\overline{\Delta_g}$ denotes the closure of
$\Delta_g$ in $L^2(X)$. Hence, the
Fourier transform yields the spectral decomposition of the selfadjoint
operator $\overline{\Delta_g}$.

\subsection{The resolvent of the Laplacian}

Let $\Hi$ be a Hilbert space and let $A$ be a closed operator in $\Hi$.  The
resolvent set of $A$ is defined to be the set of all $z \in \cz$ such
that the operator $A-z$ has a bounded inverse.  The resolvent
$R_z(A):=(A-z)^{-1}$ is a holomorphic function with values in
$\mathcal{B}(\Hi)$ on the resolvent set $\rho(A)$ of $A$. The complement of
the resolvent set is called spectrum of $A$. If $A$ is selfadjoint the
spectrum is contained in the real line.

We want to investigate the resolvent of the closure of the Laplace operator
$\Delta_g$ on a symmetric space $X=G/K$ of noncompact type with invariant
Riemannian metric $g$.  It follows from (\ref{multi}) and the above that 
the spectrum of
$\overline \Delta_g$ is exactly the set
\begin{gather}
  \spec(\overline \Delta_g)=[|\rho|^2,\infty).
\end{gather}
For $f,g \in L^2(X)$ we have by the spectral theorem
\begin{gather}
  \langle f, (\Delta_g-z)^{-1} g \rangle=
  \int_{\liea^*_+ \times B}  \frac{\overline{\hat
      f(\lambda,b)}\hat g(\lambda,b)} { |\lambda|^2 + |\rho|^2 -z}\frac{d
  \lambda}{|\hc(\lambda)|^2} db= \nonumber\\ 
  =w^{-1} \int_{\liea^* \times B} \frac{\overline{\hat
      f(\lambda,b)}\hat g(\lambda,b)} { |\lambda|^2 + |\rho|^2 -z}  \frac{d
  \lambda}{|\hc(\lambda)|^2} db\;.
\end{gather}

\section{Analytic continuation of the resolvent kernel}

Let $X=G/K$ by a symmetric space of noncompact type with invariant Riemannian
metric $g$. Then the resolvent of the Laplacian
\begin{gather}
  R_z(\overline \Delta_g)=(\overline \Delta_g - z)^{-1}
\end{gather}
can be regarded as a map $C^\infty_0(X) \to \mathcal{D}'(X)$ and hence has a
distributional kernel. This resolvent kernel is a holomorphic function on the
set $\cz \backslash [|\rho|^2,\infty)$ with values in $\mathcal{D}'(X \times
X)$.  Throughout this section $r$ will denote the positive real number
\begin{gather} \label{defr}
  r:= \mathrm{min} \{|\alpha| j(m_\alpha) \}_{\alpha \in \Delta(\lieg,\liea)_0^+},
\end{gather}
where $j(x)=\frac{x}{2}$ if $x$ is odd and $j(x)=\frac{x}{2}+1$ if $x$ is even.
Note that here $|\alpha|$ is the length of $\alpha$ with respect to the metric
$g$ which is not necessarily induced by the Killing form.  Our main result is
that the resolvent kernel admits an analytic continuation to a larger
Riemann surface.

\begin{theorem} \label{th1}
  Let $X=G/K$ be a symmetric space of noncompact type with invariant
  Riemannian metric $g$. Suppose that the rank of $X$ is odd.  Let $H^-:=\{z
  \in \cz;\; \Im(z) < 0\}$ be the lower half plane.  For each $f,g \in
  C^\infty_0(X)$ define the holomorphic function
 \begin{gather*}
   F: H^- \to \cz,\quad z \to \langle f,\; (\Delta_g
     -|\rho|^2 - z^2)^{-1} \;g \rangle.
 \end{gather*}
 Then the function $F$ admits a holomorphic continuation
 to the open set $\cz \backslash \I [r,\infty)$.
 If $\mathrm{rk}(X)=1$ the function $F$ admits a meromorphic
 continuation to $\cz$ with all poles contained in the set
 $\I [r,\infty)$.
\end{theorem}

Denote by $\Lambda^{o}$ the concrete Riemann surface associated with the function
$\sqrt{z-|\rho|^2}$ on the domain $\cz \backslash{[|\rho|^2,\infty)}$. This
means that $\Lambda^{o}$ is a branched double cover of the complex plane
with branching point $|\rho|^2$. 
The original domain $\cz \backslash{[|\rho|^2,\infty)}$ is commonly referred
to as the physical sheet, whereas its complement with the spectrum $[|\rho|^2,\infty)$
removed is called the unphysical sheet.
Let
$\Lambda^o$ be the surface with the half line $(-\infty,|\rho|^2-r^2]$
removed on the unphysical sheet. Then the above means that
the resolvent kernel regarded as a distribution extends to a function
which is holomorphic on $\Lambda_r^o$. 
In case the rank is $1$ the resolvent kernel extends to a meromorphic 
function on $\Lambda^o$ with all poles contained in the half 
line $(-\infty,|\rho|^2-r^2]$ on the unphysical sheet.
Note that in the literature sometimes another parameterization is chosen.
By a change of variables $z=\I(s-|\rho|)$ one can see that the functions
\begin{gather}
 \langle f,\left (\Delta_g - s(4 |\rho|^2-s)\right)^{-1} g \rangle,
\end{gather}
which are defined in the half plane $\Re(s) < |\rho|$ have an analytic
continuation across the line $\Re(s) = |\rho|$. 

In case the symmetric space has even rank there is an analogous result.

\begin{theorem} \label{th2}
  Let $X=G/K$ by a symmetric space of noncompact type with invariant
  Riemannian metric $g$. Suppose that the rank of $X$ is even.  Denote by
  $S_{a,b}$ the strip $\{z \in \cz;\; a < \Im(z) < b\}$.  Then for each $f,g
  \in C^\infty_0(X)$ the holomorphic function
 \begin{gather*}
   F: S_{-\pi,0} \to \cz,\quad z \to \langle f,\; (\Delta_g
     -|\rho|^2 - e^{2z})^{-1} \;g \rangle
 \end{gather*}
 admits a holomorphic continuation to the open set
 \begin{gather*}
   \mathcal{U}:=\{z \in \cz;\; z \notin \I \pi (n+\frac{1}{2}) + [\log(r),\infty)
   \quad \forall n \in (\zz \backslash \{-1\})\}.
 \end{gather*}
\end{theorem}

Let $\Lambda^e$ be the concrete Riemann surface associated with the function
$\log(z-|\rho|^2)$ on the domain $\cz \backslash{[|\rho|^2,\infty)}$.
This means, $\Lambda^e$ is the logarithmic covering of $\cz \backslash {|\rho|^2}$. 
Denote by $\Lambda_r^e$ the Riemann
surface $\Lambda^e$ with the half line $(-\infty,|\rho|^2-r^2]$
removed on all unphysical sheets. Then our result means that the resolvent
kernel regarded as a distribution has a holomorphic extension to $\Lambda_r^e$.

If all Cartan subalgebras of $G$ are conjugate there are even stronger
statements. This condition is known to be equivalent to each of the
following (see \cite{HelgI}, Ch. IX, Th. 6.1)
\begin{enumerate}
 \item all restricted roots have even multiplicity, i.e. $m_\alpha$ is even
 for all $\alpha \in \Delta(\lieg,\liea)$.
 \item $\mathrm{rk}(G)=\mathrm{rk}(X)+\mathrm{rk}(K)$.
\end{enumerate}
In this case it follows that $m_{2 \alpha}=0$ and that the function
$\hc(\lambda)^{-1}$ is a polynomial in $\lambda$ (see \cite{HelgII}, Ch. IV,
Cor. 6.15). The irreducible symmetric spaces of noncompact type, where this
happens are 
\begin{itemize}
 \item the real hyperbolic spaces of odd dimension, i.e. $SO_0(2n+1,1)/SO(2n)$,
 \item the spaces $SU^*(2n)/Sp(n)$,
 \item the spaces $G/K$ where $G$ complex,
 \item the exceptional space $E_{6(-26)}/F_4$.
\end{itemize}

\begin{theorem} \label{th3}
 Let $X=G/K$ by a symmetric space of noncompact type with invariant
 Riemannian metric $g$. Suppose that $G$ has only one conjugacy class
 of Cartan subalgebras. Let $f,g \in C^\infty_0(X)$ and
 $$
 F(z):=\langle f, (\Delta_g
     -|\rho|^2 - z)^{-1} g \rangle
 $$
 Then the following holds.
 \begin{enumerate}
  \item if  $\mathrm{rk}(X)$ is odd then $F(z^2)$
  has an analytic continuation to the whole complex plane.
  \item if  $\mathrm{rk}(X)$ is even then $F(e^{2z})$ has an analytic
  continuation to the whole complex plane.
 \end{enumerate}
\end{theorem}

Hence, in this special situation the resolvent kernel has a holomorphic extension
to $\Lambda^e$ of $\Lambda^o$, depending on whether the rank of $X$ is even or
odd. 

\section{Proof of the main results} 

We will split the proof of this theorem into several propositions.
By the spectral theorem we have for $f,g \in C^\infty_0(X)$ and $z \notin
[|\rho|^2,\infty)$
\begin{gather}
  \langle f, (\Delta_g-z)^{-1} g \rangle=  \nonumber \\
  =w^{-1} \int_{\liea^*} \frac{1}{|\hc(\lambda)|^2(|\lambda|^2 +
    |\rho|^2 -z)} \int_B
  \overline{\hat f(\lambda,b)}\hat g(\lambda,b) db \;d\lambda=\\
  = \int_{\liea^*} \frac{1}{|\hc(\lambda)|^2(|\lambda|^2 + |\rho|^2 -z)}
  V(\lambda) d \lambda, \nonumber
\end{gather}
where $V(\lambda)=w^{-1} \int_B \hat {\bar{f}} (-\lambda,b)\hat
g(\lambda,b) db$ is rapidly decaying and admits a continuation to an entire
function on $\cz$. This follows from the analytic properties of the Fourier
transforms and the fact that $B$ is compact.  Now denote by $S$ the unit
sphere in $\liea^*$, i.e. $S=\{\lambda \in \liea^*;\; |\lambda|=1\}$. Then using
polar coordinates we obtain
\begin{gather} \label{specform}
  \langle f, (\Delta_g-|\rho|^2-z)^{-1} g \rangle = \int_{\rz^+}
  \frac{F(x)}{(x^2 -z )} dx,
\end{gather}
where
\begin{gather}\label{deff}
  F(x)=C x^{\mathrm{dim}(\liea-1)} \int_S
  \frac{V(x\lambda)}{\hc(x \lambda)\hc(- x \lambda)} d \mu_S(\lambda),
\end{gather}
$\mu_S$ is the usual measure on the sphere and $C$ is a constant not depending
on $x$.

\begin{figure}
  \centerline{\includegraphics*[width=7cm]{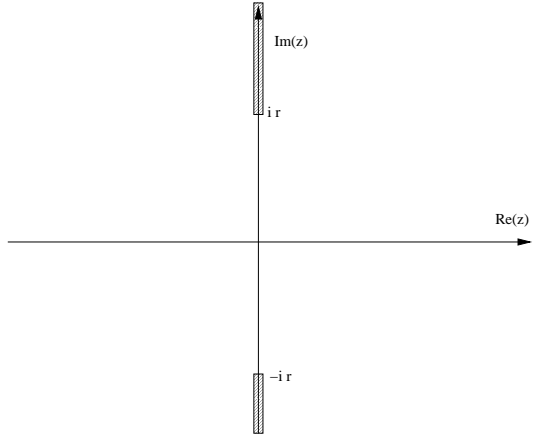}}
  \caption{Domain of analyticity of the function F(z)}\label{domain}
 \end{figure}

\begin{pro} \label{holf}
  The function $F: \rz^+ \to \cz$ defined by (\ref{deff}) is bounded and admits a
  holomorphic continuation to the set $\cz \backslash (\I [r,\infty) \cup -\I
  [r,\infty))$ (see Fig. \ref{domain}). Moreover, 
  $F(-z)=(-1)^{\mathrm{rk}(X)-1} F(z)$ and 
  $\lim_{z \to 0} z^{1-\mathrm{rk}(X)} F(z)=0$. In case $m_\alpha$ is even for all
  $\alpha \in \Delta(\lieg,\liea)_0^+$ the function $F$ admits a continuation
  to an entire function. In case the rank of $X$ is one then $F$ has a
  meromorphic extension to the whole complex plane.
\end{pro}

\begin{proof}
  The inverse of the $\hc$-function is polynomially bounded (see
  \cite{HelgII}, Ch. IV, Prop. 7.2) whereas $V$ is
  rapidly decreasing. This immediately implies that $F$ is rapidly decreasing
  as well and hence, bounded. 
  For each $\alpha \in \Delta(\lieg,\liea)_0^+$ we define the meromorphic function
  \begin{gather}
   h_\alpha(z):=
     \frac{\Gamma(\frac{1}{2}(\frac{1}{2} m_\alpha + 1 + z)) \Gamma(\frac{1}{2}(\frac{1}{2} m_\alpha +
     m_{2 \alpha} + z))}{\Gamma(z)}.
  \end{gather}
  Note that all poles of $h_\alpha$ are on the negative real axis.
  In case $m_\alpha$ is even and $m_{2 \alpha}=0$ the function $h_\alpha$
  is entire in $z$. If $m_\alpha$ is odd it is known that $m_{2
  \alpha}=0$ (see \cite{MR27:3743}, or \cite{HelgI} Chapter X, Exercise F).
  Hence, in this case the set of poles of $h_\alpha$ is
  $-(\frac{m_\alpha}{2}+n) \quad n \in \nz_0$. As in (\ref{defr}) let 
  $j(x)=\frac{x}{2}$ if $x$ is odd and $j(x)=\frac{x}{2}+1$ if $x$ is even.
  Hence, if $|z|< j(m_\alpha)$ then $z$ is not a pole of $h_\alpha$. 
  Now suppose that $|\lambda| < |\alpha| j(m_\alpha)$. 
  Then we have 
  \begin{gather}
    |\langle \lambda , \frac{\alpha}{\langle \alpha,\alpha
   \rangle} \rangle| \leq |\lambda|_\killing |\alpha|_\killing^{-1}=|\lambda|
   |\alpha|^{-1} < j(m_\alpha),
  \end{gather}
  and therefore the function $h_\alpha(\langle \I \lambda,\alpha_0 \rangle)$
  is analytic in the ball $|\lambda| < r$, where $r$ is defined by (\ref{defr}).
  By the product formula (\ref{GKFormula}) the function
  $\frac{1}{\hc(-\lambda) \hc(\lambda)}$ is analytic in the ball
  $|\lambda|<r$.
  Now suppose that $z \in \cz \backslash \I \rz$ and $\lambda \in S$. Then we have either
  $\langle \I z \lambda ,\alpha_0 \rangle=0$ or 
  $\Im (\langle \I z \lambda ,\alpha_0 \rangle) \not= 0$. In this case $z
  \lambda$ is not a pole of $h_\alpha(\langle \I \lambda,\alpha_0 \rangle)$.
  We conclude that the function $\cz \backslash (\I [r,\infty) \cup -\I
  [r,\infty)) \times S \to \cz$ 
  \begin{gather*}
   (z,\lambda) \to \frac{V(z\lambda)}{\hc(z \lambda)\hc(- z \lambda)}
  \end{gather*} 
  is holomorphic in the first variable.
  Hence, the defining formula (\ref{deff}) for
  the function $F$ yields the desired analytic continuation. The formula for
  $F(-z)$ is immediate from (\ref{deff}). Since the $\Gamma$-function has a
  pole at $0$ the function $\hc(\lambda)^{-1}$ vanishes at $0$. Hence, the
  integral on the right hand side of (\ref{deff}) vanishes at
  $\lambda=0$. This implies $\lim_{z \to 0} z^{1-\mathrm{rk}(X)} F(z)=0$.
  If all
  $m_\alpha$ are even the function $\frac{1}{\hc(-\lambda) \hc(\lambda)}$ is a
  polynomial and hence entire. In this case $F$ is entire as well by the same
  argument. In case the rank of $X$ is one the function $F$ is meromorphic
  in the whole complex plane, since in this case the function 
  $\frac{1}{\hc(-\lambda) \hc(\lambda)}$ is a meromorphic function of one
  variable and the integral in (\ref{deff}) is a sum.
\end{proof}

 By the above observation we are interested in evaluating integrals of the
 form
 $$
 \int_{\rz^+} \frac{f(x)}{x^2-z} dx,
 $$
 where $f$ is a meromorphic function on $\cz \backslash (\I [r,\infty)\cup
 -\I [r,\infty))$ with no poles on the real axis.  This integral clearly
 defines a holomorphic function on $\cz \backslash \rz^+$, and we may ask,
 whether it has an extension to a meromorphic function on a larger Riemann
 surface.  We have the following results.
\begin{pro} \label{pro1}
  Let $r$ be a positive real number and suppose that $f$ is a meromorphic
  function on the set $\cz \backslash (\I [r,\infty)\cup -\I [r,\infty))$ such
  that
 \begin{gather*}
   |f(x)| \leq C, \quad \forall x \in \rz \quad \textrm{and}\\
   f(x)=f(-x),\quad \forall x \in \rz.
 \end{gather*}
 Define the holomorphic function $G$ on $H^-:=\{z \in \cz;\; \Im(z)<0\}$ by
 $$
 G(z):= \int_{\rz} \frac{f(x)}{x^2-z^2} dx.
 $$
 Then $G(z)$ has a meromorphic continuation to $\cz \backslash \I
 [r,\infty)$. Except for the point $0$ all poles of $G(z)$ are contained in
 the set of poles of $f$. The singular behaviour of $G(z)$ at $0$ is like
 $\frac{\I \pi f(z)}{z}$. Hence, if $f(0)=0$ then $0$ is not a pole.
\end{pro}
\begin{proof}
  The function $G$ is clearly holomorphic in the lower half plane, since the
  integral converges absolutely. Let $\gamma$ be a path in $\cz$ chosen
  like in Fig. \ref{patha}, such that $\gamma$ does not meet any poles.
 \begin{figure}
   \centerline{\includegraphics*[width=7cm]{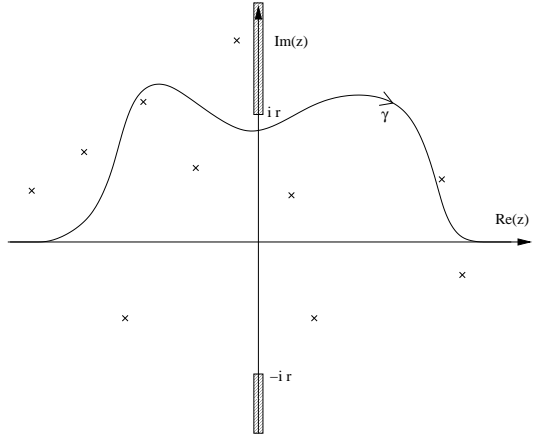}}
  \caption{The path $\gamma$, poles of $f(z)$ are indicated by crosses}\label{patha}
 \end{figure}
 Denote the domain $\cz \backslash (\I [r,\infty)\cup -\I [r,\infty))$ by
 $\mathcal{U}$.  The path divides $\mathcal{U}$ into two components. Let
 $\mathcal{U}^-_\gamma$ be the interiour of the lower component. Then
 \begin{gather}
 G_\gamma(z):= \int_{\gamma} \frac{f(x)}{x^2-z^2} dx
 \end{gather}
 defines a holomorphic function on the open set
 $O_\gamma=\mathcal{U}^-_\gamma \cap -\mathcal{U}^-_\gamma$.  Let $H^+$ and
 $H^-$ be the upper and lower half plane respectively, i.e.  $H^\pm=\{z \in
 \cz;\; \Im(z)\gtrless 0\}$.
 
 By Cauchies theorem we have for $z \in O_\gamma \cap H^-$
 \begin{gather}
 G(z)-G_\gamma(z)=2 \pi \I \sum_a \mathrm{Res}_{x=a}(\frac{f(x)}{x^2-z^2}),
 \end{gather}
 where the sum is taken over all poles $a$ contained in $O_\gamma \cap H^+$
 of the function $x \to \frac{f(x)}{x^2-z^2}$.  Now choose $z \in O_\gamma \cap H^-$
 such that $-z$ is not a pole of $f$.  Then the function
 $\frac{f(x)}{x^2-z^2}$ has a simple pole at $x=-z$ with residuum
 $\frac{f(z)}{2z}$. If $a$ is a pole of $f$ of order $k$ we have
 \begin{gather}
   \mathrm{Res}_{x=a}(\frac{f(x)}{x^2-z^2})= \frac{1}{(k-1)!}
   (\frac{\partial}{\partial x})^{k-1} \frac{f(x)}{x^2-z^2}|_{x=a}
 \end{gather}
 This expression is a finite sum of terms of the form $b_n (z^2-a^2)^{-n}$
 where the constants $b_n$ do not depend on $z$.  Hence, if $a$ is a pole of
 $f$ then $R_a(z)=\mathrm{Res}_{x=a}(\frac{f(x)}{x^2-z^2})$ is a meromorphic
 function of $z$ in the whole complex plane with poles only at $z=\pm a$. Let
 $\mathcal{P}$ be the set of poles of $f$ in $O_\gamma \cap H^+$.  We conclude
 that
 \begin{gather}
 R(z)=\frac{f(z)}{2 z} + \sum_{a \in \mathcal{P}} R_a(z)
 \end{gather}
 extends to a meromorphic function on the complex plane with poles only at
 $0$ and at the poles of $f$. The residuum of $R(z)$ at $0$ is $f(0)/2$.
 Moreover, if $z \in O_\gamma \cap H^-$ and $-z$ is not a pole of $f$ we have
 the equation
 \begin{gather}
 G(z)-G_\gamma(z)=2 \pi \I R(z).
 \end{gather}
 Since the set of such points is open this equation holds everywhere in
 $O_\gamma \cap H^-$.  Since $G_\gamma(z)$ is holomorphic in $O_\gamma$ the
 equation
 \begin{gather}
 G(z)=2 \pi \I R(z) + G_\gamma(z)
 \end{gather}
 defines a meromorphic continuation of $G(z)$ to $O_\gamma \cup H^-$.
 Since $\gamma$ can be chosen such that an arbitrary point in $\cz \backslash
 \I [r,\infty)$ is contained in $O_\gamma$ this completes the proof.
\end{proof}

This allows us to analytically continue integrals of the form
$$
\int_{\rz^+} \frac{f(x)}{x^2-z^2} dx
$$
in case $f(x)=f(-x)$. The case $f(x)=-f(-x)$ is covered by the following
proposition.

\begin{pro} \label{pro2}
  Suppose that $f: \cz \backslash (\I
  [r,\infty)\cup -\I [r,\infty)) \to \cz$ is a holomorphic function such that
 \begin{gather*}
   |f(x)| \leq C, \quad \forall x \in \rz \quad \textrm{and},\\
   f(x)=-f(-x), \quad \forall x \in \rz.
 \end{gather*}
 Let $S_{a,b}$ be the strip $\{z \in \cz;\; a < \Im(z) < b\}$ and define the
 holomorphic function $G$ on $S_{-\pi,0}$ by
 $$
 G(z):= \int_{\rz^+} \frac{f(x)}{x^2-e^{2z}} dx.
 $$
 Denote by $\mathcal{U}$ the open set
 \begin{gather*}
   \mathcal{U}:=\{z \in \cz;\; z \notin \I \pi (n+\frac{1}{2}) + [\log(r),\infty)
   \quad \forall n \in (\zz \backslash \{-1\})\}.
 \end{gather*}
 Then $G$ has an analytic continuation to $\mathcal{U}$. Moreover, denoting
 the analytic continuation again by $G$ we have
 \begin{gather*}
   G(z+\I \pi)-G(z)= \pi \I f(e^z) e^{-z}.
 \end{gather*}
\end{pro}
\begin{proof}
  A simple change of variables gives
  \begin{gather}
  G(z)=\int_{\rz} \frac{f(e^{x})}{e^{2x}-e^{2z}} e^{x} dx.
  \end{gather}
  We define $\mathcal{U}_0$ by
  \begin{gather*}
   \mathcal{U}_0:=\{z \in \cz;\; z \notin \I \pi (n+\frac{1}{2}) + [\log(r),\infty)
   \quad \forall n \in \zz\}.
 \end{gather*}
  By assumption the function $g(x):=f(e^{x}) e^{x}$ is holomorphic in
  $\mathcal{U}_0$ and is periodic in the sense that $g(x + \pi \I)=g(x)$.  Now
  choose a path $\gamma$ like in Fig. \ref{pathb}.
 \begin{figure}
   \centerline{\includegraphics*[width=7cm]{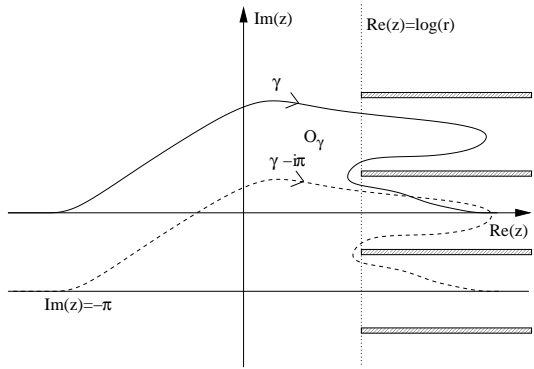}}
  \caption{The Path $\gamma$}\label{pathb}
 \end{figure}
 Denote by $O_\gamma$ the connected component of $-\infty-\I \frac{\pi}{2}$ in
 $\mathcal{U}_0 \backslash (\gamma \cup (\gamma-\I \pi \zz))$. Then
 \begin{gather}
 G_\gamma(z)=\int_{\gamma} \frac{g(x)}{e^{2x}-e^{2z}} dx.
 \end{gather}
 defines a function holomorphic in $O_\gamma$.  Now let $z$ be a point in
 $O_\gamma \cap S_{-\pi,0}$. Then by Cauchies integral theorem we have
 \begin{gather}\label{equ40}
 G(z)=G_\gamma(z).
 \end{gather}
 For each point $x
 \in \mathcal{U}_0 \cap H^+$ there is a curve, such that $x \in
 O_\gamma$. Hence, Equ. (\ref{equ40}) defines
 an analytic continuation of $G$ to $\mathcal{U} \cap H^+$. 
 Suppose now that $\gamma_1$ and $\gamma_2$ are two different paths like in 
 Fig. \ref{pathb} such that $z \in O_{\gamma_1}$ and $z+\I \pi \in
 O_{\gamma_2}$. Then Cauchies integral theorem gives
 \begin{gather}
  G_{\gamma_2}(z+\I \pi)-G_{\gamma_1}(z)=2 \pi \I \frac{g(z)}{2} e^{-2z}= \pi
  \I f(e^z) e^{-z}.
 \end{gather}
 Hence, the analytically continued function satisfies
 \begin{gather} \label{period}
  G(z+\I \pi)-G(z)= \pi \I f(e^z) e^{-z}.
 \end{gather}
 This formula provides an extension of $G$ to all of $\mathcal{U}$ satisfying
 (\ref{period}).
\end{proof}

In view of the formula (\ref{specform}) the combination of 
Prop. \ref{holf} with Prop. \ref{pro1} and Prop. \ref{pro2} proves the
Theorems \ref{th1}, \ref{th2} and \ref{th3}.

\section{Concluding Remarks}

By elliptic regularity the resolvent kernel (and its analytic continuation)
is smooth off the diagonal. Hence, we automatically obtain an analytic
continuation of the functions $R_z(x,y)$ for $x \not= y$ to the same 
Riemann surfaces.

Our result may be used to obtain meromorphic continuations of resolvent
kernels associated to operators of the form $\Delta_g + V$, where
$V$ is a compactly supported potential. This can be obtained using standard
perturbation arguments based on the analytic Fredholm theorem.

We would also like to mention that our method 
does not rely on the fact that the Fourier transforms of the
test functions are in the space of rapidly decaying functions of 
uniform exponential type $\Hi(a^*_\cz)$. 
We only use the fact that the Fourier transforms 
extend to entire function on $a^*_\cz$ and that the function $F$ defined
by $(\ref{deff})$ is integrable. The latter happens to be the case for
all $f,g \in L^2(X)$. The former requires certain growth conditions
for $f$ and $g$ at infinity. 
Suppose for example that $\alpha$ is a smooth function on $X$
such that $\alpha>1$ and such that for all $b \in B$
the functions $\alpha(x)^{-\frac{1}{2}} e^{\lambda(A(x,b))}$
are entire functions of $\lambda$ with values in $L^2(X)$.
Then our results remain true for $f,g$ in the weighted Hilbert space
 $L^2(X,\alpha(x) dx)$.
Denoting $\Hi_+=L^2(X,\alpha(x) dx)$ and $\Hi_-=L^2(X,\alpha^{-1}(x) dx)$
we therefore see that the resolvent regarded as function with values
in $\mathcal{L}(\Hi_+,\Hi_-)$ has an analytic (meromorphic) continuation
to the considered Riemann surfaces. This may be useful when considering
perturbations of the Laplace operator by potentials that are not compactly
supported.

${}$\\[10pt]
{\bf Acknowledgments}\\[6pt]

I would like to thank W. M\"uller for raising my interest in this
topic and for many discussions and hints. For the latter
I would also like to thank D. Sch\"uth and H. Thaler.

\vspace{0.5cm}

\parskip0ex

Mathematisches Institut

Universit\"at Bonn

Beringstr. 1

53115 Bonn

Germany

\vspace{0.4cm}

E-Mail: {\tt strohmai@math.uni-bonn.de}

\end{document}